\documentclass[a4paper]{amsart}

\pdfoutput=1

\usepackage{lmodern}
\usepackage{bbm}
\usepackage{amsmath,amssymb,amsthm,url,scalerel}
\usepackage[utf8]{inputenc}
\usepackage[T1]{fontenc}
\usepackage{libertine}
\usepackage[libertine,cmintegrals,cmbraces]{newtxmath}
\usepackage{verbatim}
\usepackage[shortlabels]{enumitem}
\usepackage{stmaryrd}
\usepackage{mathtools}
\usepackage{microtype}
\usepackage{multicol}
\usepackage{xcolor}
\definecolor{darkred}{RGB}{139,0,0}
\definecolor{darkblue}{RGB}{0,0,139}
\definecolor{darkgreen}{RGB}{0,100,0}
\usepackage{hyperref}
\hypersetup{
  colorlinks   = true, 
  urlcolor     = darkred, 
  linkcolor    = darkred, 
  citecolor   = darkgreen}
\usepackage{tikz}
\usepackage{tikz-cd}
\usepackage{sseq}
\usepackage[enableskew]{youngtab}
\usepackage{ytableau}
\usepackage[capitalise]{cleveref}
\usepackage{nicefrac}

\AtBeginDocument{%
   \def\MR#1{}
}
\setenumerate[1]{label=(\roman*),} 
\setitemize[1]{label=\raisebox{0.25ex}{\tiny$\bullet$}}

\newcommand{\BlockDiff}{\ensuremath{\mathrm{\widetilde{Diff}}}}

\newcommand{\BlockHomeo}{\ensuremath{\mathrm{\widetilde{Homeo}}}}

\newcommand{\hAut}{\ensuremath{\mathrm{hAut}}}

\newcommand{\OO}{\mathrm{O}}
\newcommand{\STOP}{\mathrm{STop}}

\newcommand{\id}{\mathrm{id}}

\newcommand{\G}{\mathrm{G}}

\newcommand{\BTOP}{\mathrm{BTop}}

\newcommand{\BO}{\mathrm{BO}}

\newcommand{\SO}{\mathrm{SO}}

\DeclareMathAlphabet{\mathpzc}{OT1}{pzc}{m}{it}

\newcommand{\oH}{\ensuremath{\mathrm{H}}}
\newcommand{\oG}{\ensuremath{\mathrm{G}}}
\newcommand{\oO}{\ensuremath{\mathrm{O}}}

\newcommand{\bfZ}{\ensuremath{\mathbf{Z}}}

\newcommand{\triv}{{\mathrm{triv}}}

\newcommand{\ra}{\rightarrow}

\newcommand{\xlra}[1]{\overset{#1}{\longrightarrow}}

\newcommand{\circled}[1]{\raisebox{.5pt}{\textcircled{\raisebox{-.9pt} {#1}}}}
\makeatletter
\renewcommand{\boxed}[1]{\text{\fboxsep=.2em\fbox{\m@th$\displaystyle#1$}}}
\makeatother

\newtheorem{bigthm}{Theorem}

\newtheorem*{nquestion}{Question}

\theoremstyle{definition}

\theoremstyle{remark}

\newtheorem*{nex}{Example}

\newtheorem*{nrem}{Remark}

\calclayout

\begin{document}

\title[A note on homotopy and pseudoisotopy of diffeomorphisms of $4$-manifolds]{A note on homotopy and pseudoisotopy \\of diffeomorphisms of $4$-manifolds}

\author{Manuel Krannich}
\address{Department of Mathematics, Karlsruhe Institute of Technology, 76131 Karlsruhe, Germany}
\email{krannich@kit.edu}

\author{Alexander Kupers}
\address{Department of Computer and Mathematical Sciences, University of Toronto Scarborough, 1265 Military Trail, Toronto, ON M1C 1A4, Canada}
\email{a.kupers@utoronto.ca}

\begin{abstract}This note serves to record examples of diffeomorphisms of  closed smooth $4$-manifolds $X$ that are homotopic but not pseudoisotopic to the identity, and to explain why there are no such examples when $X$ is orientable and its fundamental group is a free group.
\end{abstract}

\maketitle

Recall that two diffeomorphisms $\varphi_0$ and $\varphi_1$ of a smooth manifold $X$ are called \emph{pseudoisotopic} if there exists a diffeomorphism of $X\times [0,1]$ that restricts to $\varphi_i$ on $X\times\{i\}$ for $i=0,1$. 

\medskip

\noindent As part of the K3 project\footnote{See \url{https://aimath.org/workshops/upcoming/kirbylist/}.}, we were asked the following question:

\begin{nquestion}Is any diffeomorphism of a connected closed smooth $4$-manifold $X$ that is homotopic to the identity also pseudoisotopic to the identity?
\end{nquestion}

Equivalently, the question is whether homotopy implies pseudoisotopy for diffeomorphisms of connected closed smooth $4$-manifolds $X$. If the fundamental group $\pi_1(X)$ vanishes or if $X$ is orientable and $\pi_1(X)\cong \bfZ$, the answer to this question and its analogue in the topological category is known to be positive (see \cite[Theorem 1]{Kreck}, \cite[Proposition 2.2]{Quinn}, and  \cite[p.\,51]{StongWang}).

\medskip

\noindent The purpose of this note is twofold: firstly, we illustrate how classical surgery theory allows one to answer the general form of the question in the negative. 

\begin{bigthm}\label{thm:non-example}
	There exists a diffeomorphism of a smooth closed $4$-manifold $X$ that is homotopic but neither smoothly nor topologically pseudoisotopic to the identity.
\end{bigthm}

More concretely, we explain why a diffeomorphism as in \cref{thm:non-example} exists for any $4$-manifold of the form $X=Y\sharp^g(S^2\times S^2)$ for large enough $g\ge0$ depending on $Y$, where $Y$ is any connected compact smooth stably parallelisable $4$-manifold whose fundamental group $\pi\coloneq \pi_1(Y)$ satisfies:
	\begin{enumerate}[(1)]
		\item\label{enum:condition-pi-1} $\oH_1(\pi)/(\mathrm{2}\text{-torsion})$ is not annihilated by multiplication by $3$ and
		\item\label{enum:condition-pi-2} the $5$th simple $L$-group $L_5^s(\bfZ[\pi])$ of $\bfZ[\pi]$ with the standard involution vanishes.
	\end{enumerate}
There are many $4$-manifolds with these properties: any finitely presented group $\pi$ arises as $\pi_1(Y)$ of a connected closed smooth stably parallelisable $4$-manifold $Y$ \cite[Proof of Theorem 1]{Kervaire} and there are many choices for $\pi$ that satisfy \ref{enum:condition-pi-2}, such as finite groups $\pi$ of odd order \cite{Bak} or more generally products $\pi=\pi_{\mathrm{odd}}\times\bfZ/2^k$ where $\pi_{\mathrm{odd}}$ has odd order \cite[p.~227, 12.1, 12.2]{HambletonTaylor}, e.g.\,$\pi$ can be any finite cyclic group. If $\pi_{\mathrm{odd}}$ has a nontrivial element of order $\neq 3$, then  $\pi$ also satisfies \ref{enum:condition-pi-1}.

\begin{nex}The simplest example for $Y$ that satisfies the conditions is the result of a surgery along an embedding $e\colon S^1\times D^3\hookrightarrow S^1\times S^3$ such that (a) the class $[e]\in\pi_1(S^1\times S^3)\cong\bfZ$ is $\pm5$ and (b) the result of the surgery is stably parallelisable (which is always possible; see \cite[Theorem 2]{MilnorHomotopy}).
\end{nex}

The second purpose of this note is to observe that a combination of the surgery exact sequence with arguments in work of Shaneson allows one to widen the class of examples for which the answer to the Question is positive from orientable $4$-manifolds $X$ such that $\pi_1(X)$ is trivial or free of rank $1$ to those for which $\pi_1(X)$ is a free group $F_n$ of arbitrary finite rank $n\ge0$.

\begin{bigthm}\label{thm:fn-example} For diffeomorphisms of connected closed smooth orientable $4$-manifolds with free fundamental group, homotopy implies pseudoisotopy.
\end{bigthm}

\begin{nrem}Combined with recent work of Gabai \cite[Theorem 2.5, Remark 2.10]{GabaiSchoenflies}, this implies that diffeomorphisms of $4$-manifolds $X$ as in \cref{thm:fn-example} that are homotopic are also stably isotopic, i.e.\,are isotopic when extended by the identity to diffeomorphisms of  $X\sharp^g(S^2\times S^2)$ for some $g\ge0$ after isotoping them to fix an embedded disc to form the connected sum.
\end{nrem}

\subsection*{Proof of \cref{thm:non-example}}
The proof centres around the diagram of groups
\vspace{-0.1cm}
\begin{equation}\label{equ:big-diagram}
\begin{tikzcd}[ar symbol/.style = {draw=none,"\textstyle#1" description,sloped},
	subset/.style = {ar symbol={\subseteq}}, row sep=0.5cm,column sep=0.4cm]
\pi_0\hAut^s_\partial(X\times I)\arrow[r,"\circled{1}"]& S_\partial^{s,\triv}(X\times I)\arrow[r,"\circled{2}"]\arrow[d,subset]& \pi_0\BlockDiff_\partial(X)\arrow[r,"\circled{3}"]& \pi_0\hAut^s_\partial(X)\\[-0.3cm]
&S_\partial^{s}(X\times I)\arrow[d,"\circled{4}"]&\\
&\left[\Sigma(X/\partial X),\oG/\oO\right]_\ast\arrow[r,"\circled{6}"]\arrow[d,"\circled{5}"]& \left[\Sigma(X/\partial X),\BO\right]_\ast\\
&L^s_5(\bfZ[\pi_1(X)],w_1)&
\end{tikzcd}\vspace{-0.1cm}
\end{equation}
for any compact connected smooth $4$-manifold $X$. The terms and maps involved are:
\begin{enumerate}[leftmargin=0.6cm,label=(\Roman*) ]
\item  $\smash{\pi_0\BlockDiff_\partial(X)}$ and $\pi_0\hAut^s_\partial(X)$ are the groups of diffeomorphisms (or simple homotopy equivalences) of $X$ that fix $\partial X$ pointwise, up to smooth pseudoisotopy (or homotopy) fixing $\partial X$,
\item\label{enum:structure-set} $S_\partial^{s}(X\times I)$ is the \emph{simple structure set} of the pair $(X\times I,\partial (X\times I))$ in the sense of surgery theory \cite[Chapter 10]{Wall}, consisting of equivalence classes of pairs $(W,\varphi)$ of a compact $5$-manifold $W$ and a simple homotopy equivalence $\varphi\colon W\ra X\times I$ that restricts to a diffeomorphism $\varphi|_{\partial W}\colon \partial W\ra \partial (X\times I)$. The group structure is by ``stacking''. Note that the manifold $W$ becomes via $\varphi|_{\partial W}$ a relative self $s$-cobordism of $(X,\partial X)$ with an identification of the boundary bordism with $\partial X\times I$, but since $W$ is $5$-dimensional this $s$-cobordism need not be trivial. Restricting to those classes of pairs $(W,\varphi)$ for which $W$ \emph{is} trivial (i.e.\,diffeomorphic to $X\times I$ relative to $\partial X\times I$, but not necessarily relative to $X\times \{0,1\}$) defines a subgroup $S_\partial^{s,\triv}(X\times I)\subseteq S_\partial^{s}(X\times I) $.
\item $[\Sigma(X/\partial X),\BO]_\ast$ and $[\Sigma(X/\partial X),\oG/\oO]_\ast$ consist of pointed homotopy classes of maps from the reduced suspension of $X/\partial X$ to the classifying spaces $\BO$ and $\oG/\oO$ for stable vector bundles and for stable vector bundles with a trivialisation of the underlying stable spherical fibration,
\item $L^s_5(\bfZ[\pi_1(X)],w_1)$ is the $5$th simple $L$-group of the group ring $\bfZ[\pi_1(X)]$ with the involution determined by orientation character $w_1$, in the sense of surgery theory (see loc.cit.).
\item\label{enum:maps} The map $\smash{\circled{1}}$ sends $[\varphi]$ to $[X\times I, \varphi]$, $\smash{\circled{2}}$ sends $[X\times I,\varphi]$ to $\smash{[\varphi|_{X\times \{1\}}\circ \varphi|_{X\times \{0\}}^{-1}]}$, $\smash{\circled{3}}$ sends $[\varphi]$ to $[\varphi]$, $\smash{\circled{4}}$ and $\smash{\circled{5}}$ take normal invariants and surgery obstructions respectively (see loc.cit.), and $\smash{\circled{6}}$ is induced by forgetting the trivialisation of the spherical fibration.
\end{enumerate}

\noindent We will use the following facts about this diagram:
\begin{enumerate}[leftmargin=0.6cm,label=(\alph*)]
\item\label{enum:exact} The top row and middle column are exact; the former by inspection and the latter by surgery theory (see e.g.\,Chapter 10 loc.cit.).
\item \label{enum:normal-invariant}  We have $[\Sigma(X/\partial X),\oG/\oO]_\ast\cong \oH^1(X,\partial X;\bfZ/2)\oplus  \oH^3(X,\partial X;\bfZ)$ (see e.g.\,\cite[p.\,398]{KirbyTaylor}). 
\item\label{enum:geometric-map} The composition $S_\partial^{s}(X\times I)\ra [\Sigma(X/\partial X),\BO]_\ast$ has the following description (cf.\,\cite[p.\,113--114]{Wall}): for $[W,\varphi]\in S_\partial^{s}(X\times I)$, choose a homotopy inverse $\widetilde{\varphi}$ of $\varphi$ with $\widetilde{\varphi}|_{\partial(X\times I)}=\varphi|_{\partial W}^{-1}$, consider the stable vector bundle $\widetilde{\varphi}^*(\nu_W)\oplus\tau_{X\times I}$ over $X\times I$ where $\nu_{(-)}$ and $\tau_{(-)}$ is the stable normal and tangent bundle. Together with the trivialisation of $(\widetilde{\varphi}^*(\nu_W)\oplus\tau_{X\times I})|_{\partial(X\times I)}$ induced by the derivative of $\widetilde{\varphi}_{\partial(X\times I)}$, this defines a class in $[(X\times I)/\partial (X\times I),\BO]_\ast=[\Sigma(X/\partial X),\BO]_\ast$. 

\item\label{enum:extbyid} If $X$ arises as a codimension $0$ submanifold $\smash{X\subset \mathrm{int}(\overline{X})}$ of the interior of a compact connected smooth $4$-manifold $\smash{\overline{X}}$, there is a morphism $S_\partial^s(X\times I)\ra S_\partial^s(\smash{\overline{X}}\times I)$ given by sending $(W,\varphi)$ to $W\cup_{\partial X\times I}(\smash{\overline{X}}\backslash \mathrm{int}(X)\times I)$ and  $\varphi$ to $(\varphi\cup_{\partial X\times I}\id)$. This preserves the $\triv$-subgroups from \ref{enum:structure-set}. Moreover, from \ref{enum:geometric-map} one sees that this morphism is compatible with the map $[\Sigma(X/\partial X),\BO]_\ast\ra [\Sigma(\smash{\overline{X}}/\partial \smash{\overline{X}}),\BO]_\ast$ given by precomposition with the map $\Sigma(\smash{\overline{X}}/\partial \smash{\overline{X}})\ra \Sigma(X/\partial X)$ induced by collapsing $\smash{\overline{X}}\backslash\mathrm{int}(X)$. In particular, if $\partial X\neq \varnothing$, there is a stabilisation morphism $S_\partial^{s}(X\times I)\ra S_\partial^{s}(X\sharp(S^2\times S^2)\times I)$ by viewing $X$ as a submanifold of $X\cup_{\partial X}((\partial X\times I)\sharp S^2\times S^2)\cong X\sharp(S^2\times S^2)$. 
\item \label{enum:stable-scob} If $\partial X\neq \varnothing$, then by the stable $s$-cobordism theorem \cite[1.1]{Quinnscob}, for any element $x\in  S_\partial^{s}(X\times I)$ there is a number $g\ge 0$ depending on $x$ such that the image of $x$ under the iterated stabilisation map $\smash{ S_\partial^{s}(X \times I)\ra S_\partial^{s}(X\sharp(S^2\times S^2)^{\sharp g} \times I)}$ from \ref{enum:extbyid} lies in the $\triv$-subgroup from \ref{enum:structure-set}.
\item\label{enum:localise} After localisation $\smash{(-)[\tfrac{1}{2}]}$ away from $2$ and Postnikov $7$-truncation $\tau_{\le 7}(-)$ the spaces $\oG/\oO$ and $\BO$ are both equivalent to $K(\bfZ[\tfrac{1}{2}],4)$ and the map $\tau_{\le 7}(\oG/\oO)[\tfrac{1}{2}]\ra \tau_{\le 7}(\BO)[\tfrac{1}{2}]$ is induced by multiplication by $\pm3$. This follows from the computation of the stable homotopy groups of spheres in small degrees and the surjectivity of the stable $J$-homomorphism in degree $3$.
\end{enumerate}

\noindent If $X$ is stably parallelisable and $\pi_1(X)$ satisfies the conditions \ref{enum:condition-pi-1} and \ref{enum:condition-pi-2} (let us call such manifolds \emph{admissible} for brevity), we can add to this list:
\begin{enumerate}[leftmargin=0.6cm,label=(\alph*),resume]
\item\label{enum:trivial-L} The map $\circled{4}$ is surjective as $L^s_5(\bfZ[\pi_1(X)],w_1)=0$ by \ref{enum:condition-pi-2} since $X$ is orientable.
\item \label{enum:trivial-action}The composition $(\circled{6}\circ \circled{4}\circ \circled{1})\colon \pi_0\hAut^s_\partial(X\times I)\ra [\Sigma(X/\partial X),\BO]_*$ is trivial. This follows from the descriptions of $\circled{1}$ and $\circled{6}\circ \circled{4}$ in \ref{enum:maps} and \ref{enum:geometric-map} using that $X$  is stably parallelisable.
\item\label{enum:nontrivial} By \ref{enum:normal-invariant} and Poincaré duality we have $\smash{[\Sigma(X/\partial X),\oG/\oO]_*\cong\oH_1(X;\bfZ)\oplus\oH^1(X,\partial X;\bfZ/2)}$. As a result of \ref{enum:condition-pi-1} there exists an element $x\in \oH_1(X;\bfZ)$ with $3\cdot x\neq0\in \oH_1(X;\bfZ[\tfrac{1}{2}])$. By \ref{enum:trivial-L} the class $x$ lifts along $\circled{4}$ to $\smash{S_\partial^{s}(X\times I)}$ and since  $\smash{[\Sigma(X/\partial X),\oG/\oO[\tfrac{1}{2}]]_*\ra [\Sigma(X/\partial X),\BO[\tfrac{1}{2}]]_*}$ is given by multiplication by $3$ in view of \ref{enum:localise}, it follows that $x$ maps nontrivially to $[\Sigma(X/\partial X),\BO[\tfrac{1}{2}]]_*$. In particular the composition $\smash{S_\partial^{s}(X\times I)\ra [\Sigma(X/\partial X),\BO[\tfrac{1}{2}]]_*\cong\oH_1(X;\bfZ[\tfrac{1}{2}])}$ is nontrivial.
\end{enumerate}

\noindent 
Combining all this, the proof goes as follows: given an admissible $4$-manifold $Y$, the claim is that there is a $g\ge0$ and a class in the kernel of $\smash{\circled{3}}$ for $X=Y\sharp^g(S^2\times S^2)$ that maps nontrivially to the group $\smash{\pi_0\BlockHomeo(X)}$ of topological pseudoisotopy classes of homeomorphisms. We fix an embedded disc $D^4\subset \mathrm{int}(Y)$ and set $Y^{\circ}\coloneq Y\backslash\mathrm{int}(D^4)$ which is again admissible. By \ref{enum:nontrivial} there is a class $x\in\smash{S_\partial^{s}(Y^\circ \times I)}$ that maps nontrivially to  $[\Sigma(Y^\circ/\partial Y^\circ),\BO[\tfrac{1}{2}]]_*$. Using \ref{enum:stable-scob}, since $\partial Y^\circ \neq \varnothing$ we find $g\ge0$ such that the image $x_g$ of $x$ under $g$-fold stabilisation lies in $\smash{S_\partial^{s,\triv}(Y^\circ\sharp(S^2\times S^2)^{\sharp g} \times I)}$. Using the inclusion $\smash{Y^\circ\sharp(S^2\times S^2)^{\sharp g}\subset Y\sharp(S^2\times S^2)^{\sharp g}}$ we obtain via \ref{enum:extbyid} an element $\overline{x_g}$ in $\smash{S_\partial^{s,\triv}(Y\sharp^g(S^2\times S^2) \times I)}$. By the compatibility part of \ref{enum:extbyid} the image of this element in $\smash{[\Sigma(Y\sharp(S^2\times S^2)^{\sharp g}/\partial(Y\sharp(S^2\times S^2)^{\sharp g})),\BO[\tfrac{1}{2}]]_*}$ is the image of $x$ under the composition $\smash{S_\partial^{s}(Y^{\circ}\times I) \ra [\Sigma(Y^{\circ}/\partial Y^{\circ}),\BO[\tfrac{1}{2}]]_*\ra [\Sigma(Y\sharp^g(S^2\times S^2)/\partial(Y\sharp(S^2\times S^2)^{\sharp g})),\BO[\tfrac{1}{2}]]_*}$, so it is nontrivial since it is not in the kernel of the first map by choice of $x$ and because the second map is in light of \ref{enum:localise},  \ref{enum:extbyid} and Poincaré duality an isomorphism since the inclusion $Y^{\circ}\subset Y\sharp^g(S^2\times S^2)$ is an isomorphism on first homology. Using \ref{enum:trivial-action} this implies that $\overline{x_g}$ is not in the image of $\circled{1}$ for $X\coloneq Y\sharp(S^2\times S^2)^{\sharp g}$, so its image under $\circled{2}$ is a nontrivial element in the kernel of $\circled{3}$. To see that this element is also nontrivial in $\smash{\pi_0\BlockHomeo(X)}$, one argues as follows: by forgetting smoothness the diagram \eqref{equ:big-diagram} maps compatibly to the corresponding diagram in the topological category, so it suffices to show that $\overline{x_g}$ is, when mapped to the topological version of $\smash{S_{\partial}^s(X\times I)}$, still not hit by the topological analogue of $\circled{1}$. By the way we detected this element, it suffices to show that the map $[\Sigma(X/\partial X),\BO[\tfrac{1}{2}]]_\ast\ra [\Sigma(X/\partial X),\BTOP[\tfrac{1}{2}]]_\ast$ induced by the map $\BO\ra\BTOP$ classifying the underlying stable Euclidean space bundle of a stable vector bundle is injective. This holds because $\tau_{\le 7}(\BO)[\tfrac{1}{2}]\ra\tau_{\le 7}(\BTOP)[\tfrac{1}{2}]$ is an equivalence (see e.g.\,\cite[p.\,246, 5.0, (5)]{KirbySiebenmann}).
 \vspace{-0.1cm}

\begin{nrem}In dimensions $d\ge5$, it is significantly easier to produce examples as in \cref{thm:non-example}, even $1$-connected ones, for example $X=S^3\times S^n$ for any $n\ge2$ works. This can be shown by a variant of the strategy above, but for $n\ge3$ there is also a more elementary argument: choose $[\varphi]\in \pi_3(\SO(n+1))$ that maps to a nontorsion class in $\pi_3(\SO)\cong\bfZ$ and lies in the kernel of the map induced by the forgetful map $\SO(n+1)\ra \hAut(S^n)$ (this is possible if and only if $n\ge3$) and consider the orientation-preserving diffeomorphism $t_\varphi$ of $S^3\times S^n$ given by $t_{\varphi}(x,v)=(x,\varphi(x)\cdot v)$. Since the image of $\varphi$ in $\pi_3(\hAut(S^n))$ is trivial, $t_{\varphi}$ is homotopic to the identity. With respect to the standard stable framing of $S^3\times S^n$, the stable derivative of $t_{\varphi}$ in $[S^3\times S^n,\SO]\cong \pi_3(\SO) \oplus \pi_n(\SO) \oplus \pi_{n+3}(\SO)$ is given by the stabilisation of $[\varphi]$ in the first term, so it is nontrivial by choice of $\varphi$. This implies that $t_\varphi$ is not pseudoisotopic to the identity, and since the image of $\varphi$ in $\pi_3(\STOP)$ is nontrivial as $\bfZ \cong \pi_3(\SO)\ra \pi_3(\STOP)\cong \bfZ\oplus \bfZ/2$ is injective, $t_{\varphi}$ is also not topologically pseudoisotopic to the identity.
\end{nrem}

\subsection*{Proof of \cref{thm:fn-example}}
For a connected closed smooth orientable $4$-manifold $X$ with $\pi_1(X) \cong F_n$ for some $n\ge0$, the claim is that the map $\smash{\circled{3}}$ in  \eqref{equ:big-diagram} is injective which is by exactness of the top row equivalent to showing that $\smash{\circled{1}}$ is surjective. To see this, we use that---since they are part of the surgery exact sequence (see e.g.\,\cite[Chapter 10]{Wall})---the maps $\smash{\circled{4}}$ and $\smash{\circled{5}}$ can be extended to the left to an exact sequence of groups (note that since $X$ is orientable and closed, the involution on $\bfZ[\pi_1(X)]$ is the standard one and we have $X/\partial X=X_+\coloneq X\sqcup\{\ast\}$)\vspace{0.1cm}
\begin{equation}\label{equ:extended-surgery}
\smash{L_6^s(\bfZ[\pi_1(X)]) \xlra{\circled{7}} S^s_\partial(X \times I) \xlra{\circled{4}} [\Sigma(X_+),\G/\OO]_*\xlra{\circled{5}} L_5^s(\bfZ[\pi_1(X)])}.\end{equation}
The map $\circled{7}$ turns out to be trivial: choosing an embedded disc $D^4\subset X$ and using the naturality of the surgery exact sequence in codimension $0$ embeddings, this follows by combining that $L_6^s(\bfZ[1]) \to L_6^s(\bfZ[F_n])$ is an isomorphism \cite[Corollary 6]{Cappell} and that $S^s_\partial(D^4\times I)=0$ as a consequence of the solution of the $5$-dimensional smooth Poincaré conjecture. Since $\smash{\circled{7}}$ is trivial, the map $\circled{4}$ is injective. In the case $n=0$, so if $\pi_1(X)$ and thus $\oH_1(X)$ vanish, $[\Sigma(X_+),\oG/\oO]_*$ vanishes by \ref{enum:normal-invariant} and Poincar\'e duality, so it follows from exactness of \eqref{equ:extended-surgery} that $S^s_\partial(X \times I)$ vanishes, so in particular $\circled{1}$ is surjective and the proof is finished in this case.

For $n>0$, the proof that $\smash{\circled{1}}$ is surjective is more subtle and follows by adapting arguments of Shaneson in the case $n=1$: as a first step, one argues as in \cite[p.~349]{ShanesonLow} that the composition $\smash{\circled{6}\circ\circled{4}}$ is trivial: By ``gluing the ends of $X\times I$'' the composition features in a commutative diagram 
\[\begin{tikzcd}[row sep=0.5cm] S^s_\partial(X \times I) \rar{\circled{4}} \dar{\circled{8}} & {[\Sigma(X_+),\oG/\oO]_*} \dar \rar{\circled{6}} & {[\Sigma(X_+),\BO]_*} \dar \\
S^s(X \times S^1) \rar & {[(X \times S^1)_+,\oG/\oO]_*} \rar & {[(X \times S^1)_+,\BO]_*}\end{tikzcd}\]
whose middle and right vertical map are split injective, using that the quotient $(X\times S^1)_+\ra\Sigma (X_+)$ splits after suspension and that $\G/\OO$ and $\BO$ are loop spaces. By the description of $\circled{6}\circ\circled{4}$ from \ref{enum:geometric-map}, it suffices to show that any equivalence $W\ra X\times S^1$ from a closed smooth $5$-manifold preserves the stable tangent bundle. This follows from the proof of \cite[Theorem 6.1]{ShanesonProduct} (the statement assumes that the fundamental group is free abelian, but the proof goes through for $\pi_1(X\times S^1)\cong \bfZ\times F_n$ since the version of the Novikov conjecture proved in \cite[Theorem 7]{FarrellHsiang} applies).

Since $\smash{\circled{6}\circ\circled{4}}$ is trivial, it follows from \ref{enum:localise} and the fact that $\oH^3(X;\bfZ)\cong\oH_1(X;\bfZ)\cong\bfZ^n$ is torsion free that $\smash{\circled{4}}$ lands in the $2$-torsion subgroup of $ [\Sigma(X_+),\oG/\oO]_*$ which is isomorphic to $\oH^1(X;\bfZ/2)$. By exactness of the upper row in \eqref{equ:big-diagram} and injectivity of $\circled{4}$, in order to show that $\smash{\circled{1}}$ is surjective it suffices to prove that $\smash{(\circled{4}\circ\circled{1})\colon \pi_0\hAut^s_\partial(X \times I)\to [\Sigma(X_+),\oG/\oO]_*}$ surjects onto the $2$-torsion subgroup of $[\Sigma(X_+),\oG/\oO]_*$. This follows from constructing homotopy equivalences as in \cite[p.~349--350]{ShanesonLow}: first one shows that the Hurewicz homomorphism $\pi_3(X)\ra \oH_3(X) \cong \bfZ^n$ is surjective, which can be done as in the proof of Lemma 6.2 loc.cit.\,using the $1$-truncation $X \to B(F_n)\simeq \vee_n S^1$, then one chooses a basis $(\beta_i)$ for 
$\oH_3(X)$, lifts each $\beta_i$ to $\pi_3(X)$ and use the lifts to construct elements $h_i\in \pi_0\hAut^s_\partial(X \times I)$ analogous to the construction of $h$ on the top of p.~350 loc.cit. Since $\circled{4}$ is injective, if a sum $\sum_i \varepsilon_i h_i\in \pi_0\hAut^s_\partial(X \times I)$ with $\varepsilon_i \in \{0,1\}$ not all zero were in the kernel of $\circled{4}\circ\circled{1}$, then its image under $\circled{8}$ would be homotopic to a diffeomorphism. But the argument on p.~350 loc.cit.~shows that this is not the case, so the images of the $h_i$ in $[\Sigma(X_+),\oG/\oO]_*\cong \oH^1(X;\bfZ/2)$ are linearly independent and hence $\circled{4}\circ\circled{1}$ is surjective for dimension reasons.

\begin{nrem}\cref{thm:fn-example} was proved for \emph{closed} manifolds, but it does not seem unreasonable that the proof extends to allow nonempty boundary. For the boundary connected sums $\smash{Y_g\coloneq \natural^g(S^1\times D^3)}$ with $g\ge0$ the statement can in fact be \emph{deduced} from the closed case by a trick:

Given a diffeomorphism $\phi\colon Y_g\ra Y_g$ that fixes $\partial Y_g$ and is homotopic to $\id_{Y_g}$, we will show that $\phi$ is also pseudoisotopic to $\id_{Y_g}$ relative to $\partial Y_g$. Extending $\phi$ by the identity on the second copy of $Y_g$ in its double $\smash{D(Y_g)\coloneq Y_g\cup_\partial \overline{Y}_g\cong\sharp^g(S^1\times S^3)}$, we obtain a diffeomorphism $(\phi\cup \id)\colon D(Y_g)\ra D(Y_g)$ that is homotopic to the identity. As $\pi_1(D(Y_g))\cong F_g$, \cref{thm:fn-example} ensures the existence of a pseudoisotopy $H$ from $(\phi\cup \id)$ to $\id_{D(Y_g)}$. Restricting $H$ to $\smash{\overline{Y}_g}$ gives a concordance embedding $\smash{H|_{I\times \overline{Y}_g}}$ from $\smash{\overline{Y}_g}$ into $D(X)$, and by isotopy extension it suffices to show that $\smash{H|_{I\times \overline{Y}_g}}$ is isotopic, as a concordance embedding, to the inclusion $\smash{I\times \overline{Y}_g\subset I\times D(Y_g)}$. As $\smash{\overline{Y}_g\subset D(Y_g)}$ has handle codimension $\geq 3$, this holds by \cite[Theorem 2.1, Addendum 2.1.2]{Hudson}.
\end{nrem}

\subsection*{Acknowledgements}We are grateful to Oscar Randal-Williams for pointing out an oversight in an earlier version of this note. AK acknowledges the support of the Natural Sciences and Engineering Research Council of Canada (NSERC) [funding reference number 512156 and 512250]. AK was supported by an Alfred P.~Sloan Research Fellowship.

\vspace{-0.1cm}
\bibliographystyle{amsalpha}
\bibliography{literature}

\end{document}